\documentclass[12pt,centertags]{amsart}
\usepackage{amsmath,amstext,amsthm,a4,amssymb,amscd}
\usepackage[mathscr]{eucal}
\usepackage{mathrsfs}
\newcommand{\field}[1]{\mathbb{#1}}

\newcommand{\R}{\field{R}}
\newcommand{\C}{\field{C}}
\newcommand{\N}{\field{N}}

 \def\cC{\mathscr{C}}
\def\cL{\mathscr{L}}
\def\cO{\mathscr{O}}

\def\mQ{\mathcal{Q}}

\DeclareMathOperator{\End}{End}

\DeclareMathOperator{\Ker}{Ker}

\newcommand{\om}{\omega}

\newtheorem{thm}{Theorem}

\theoremstyle{definition}

\theoremstyle{definition}

\newcommand{\be}{\begin{eqnarray}}
\newcommand{\ee}{\end{eqnarray}}

\newcommand{\var}{\varepsilon}

\newcommand{\comment}[1]{}

\begin{document}

\title{A remark on 'Some numerical results in complex differential geometry'}
\date{\today}
\author{Kefeng Liu}
\address{Center of Mathematical Science, Zhejiang University and Department of Mathematics, UCLA, CA 90095-1555,
USA (liu@math.ucla.edu)}

\author{Xiaonan Ma}
\address{Centre de Math\'ematiques Laurent Schwartz, UMR 7640 du CNRS,
Ecole Polytechnique, 91128 Palaiseau Cedex,
France (ma@math.polytechnique.fr)}

\begin{abstract}
In this note we verify certain statement about the operator $Q_K$
constructed by Donaldson in \cite{D05} by using the full asymptotic
expansion of Bergman kernel obtained in \cite{DLM04a} and
\cite{MM04a}.
\end{abstract}

\maketitle

 In order to find explicit numerical approximation of
K\"ahler-Einstein metric of projective manifolds, Donaldson
introduced in \cite{D05} various operators with good properties to
approximate classical operators. See the discussions in Section
4.2 of \cite{D05} for more details related to our discussion. In
this note we verify certain statement of Donaldson about the
operator $Q_K$ in Section 4.2 by using the full asymptotic
expansion of Bergman kernel derived in \cite[Theorem 4.18]{DLM04a}
and \cite[\S 3.4]{MM04a}. Such statement is needed for the
convergence of the approximation procedure. 

Let $(X,\omega, J)$ be a compact K\"ahler manifold of
$\dim_{\C}X=n$, and let $(L, h^L)$ be a holomorphic Hermitian line
bundle on $X$. Let $\nabla ^L$ be  the holomorphic Hermitian
connection on $(L,h^L)$ with curvature $R^L$. We assume that
 \begin{align} \label{n1}
\frac{\sqrt{-1}}{2 \pi} R^L=\om.
\end{align}

Let $g^{TX}(\cdot,\cdot):= \om(\cdot,J\cdot)$ be the Riemannian metric on $TX$
induced by $\om, J$. Let $dv_X$ be the Riemannian volume
form of $(TX, g^{TX})$, then $dv_X= \om^n/n!$.
Let $d\nu$ be any volume form on $X$.
Let $\eta$ be the positive function on $X$ defined by
 \begin{align} \label{n3}
dv_X= \eta\, d\nu.
\end{align}
The $L^2$--scalar product $\langle \quad\rangle_\nu$ on $\cC^\infty (X,L^p)$,
the space of smooth sections of $L^p$, is given by
\begin{align}\label{n2}
&\langle \sigma_1,\sigma_2 \rangle_\nu
:=\int_X\langle \sigma_1(x),\sigma_2(x)\rangle_{L^p}\,d\nu(x)\,.
\end{align}

Let $P_{\nu, p}(x,x')$ $(x,x'\in X)$ be the smooth kernel of
the orthogonal projection from
$(\cC ^\infty(X, L^p),\langle\quad \rangle_\nu )$
onto $H^0(X,L^p)$, the space of the holomorphic sections of $L^p$ on $X$,
 with respect to $d\nu (x')$.
Note that $P_{\nu, p}(x,x')\in L^p_x\otimes L^{p*}_{x'}$.
Following \cite[\S 4]{D05}, set
\begin{align}\label{1n1}
K_p(x,x'):= |P_{\nu, p}(x,x')|^2_{h^{L^p}_x\otimes h^{L^{p*}}_{x'}}, \quad
R_p:= (\dim H^0(X, L^p))/ {\rm Vol}(X,\nu),
\end{align}
here ${\rm Vol}(X,\nu):= \int_X d\nu$. Set ${\rm Vol}(X,dv_X):= \int_X dv_X$.

Let $Q_{K_p}$ be the integral operator associated to $K_p$
which is defined for $f\in \cC^\infty (X)$,
\begin{equation}\label{1n2}
Q_{K_p} (f)(x):= \frac{1}{R_p} \int_X K_p(x,y) f(y)d\nu(y).
\end{equation}

Let $\Delta$ be the (positive) Laplace operator on $(X, g^{TX})$
acting on the functions on $X$.
We denote by $|\quad|_{L^2}$ the $L^2$-norm on the function
on $X$ with respect to $dv_X$.

\comment{
The following Proposition verifies that the operator
\begin{align}\label{1n3}
W(\eta)=\eta- \frac{\kappa}{2p}\Delta Q_{K_p}(\eta)
\end{align}
defined in \cite[\S 4.2]{D05}
is a contraction for $p$ large enough and $\kappa$ small.
}

\begin{thm} \label{nt1} There exists a constant $C>0$ such that for any
$f\in \cC^\infty (X)$, $p\in \N$,
\begin{align}\label{1n4}
\begin{split}
& \left|\Big( Q_{K_p}-\frac{{\rm Vol}(X, \nu)}{{\rm Vol}(X, dv_X)}
\eta \exp\Big(-\frac{\Delta}{4\pi p}\Big)\Big) f\right|_{L^2}
\leq \frac{C}{p}    \left|f\right|_{L^2},\\
&  \left| \Big(\frac{\Delta}{p} Q_{K_p}
- \frac{{\rm Vol}(X, \nu)}{{\rm Vol}(X, dv_X)}
\frac{\Delta}{p}\eta \exp\Big(-\frac{\Delta}{4\pi p}\Big)\Big) f\right|_{L^2}
\leq \frac{C}{p}    \left|f\right|_{L^2}.
\end{split}
\end{align}
Moreover, \eqref{1n4} is uniform in that there is an integer $s$
such that if all data $h^L$, $d\nu$ run over a set which are
bounded in $\cC^s$ and that $g^{TX}$, $dv_X$ are bounded from
below, then the constant $C$ is independent of $h^L$, $d\nu$.
\end{thm}
\begin{proof}
We explain at first the full asymptotic expansion of  $P_{\nu, p}(x,x')$
from \cite[Theorem 4.18$^\prime$]{DLM04a} and \cite[\S 3.4]{MM04a}.
For more details on our approach we also refer the readers
to the recent book \cite{MM06a}.

Let $E=\C$ be the trivial holomorphic line bundle on $X$. Let
$h^E$ the metric on $E$ defined by $|1|_{h^E}^2= 1$, here $1$ is
the canonical unity element of $E$. We identify canonically $L^p$
to $L^p\otimes E$ by Section $1$.

As in \cite[\S 3.4]{MM04a},
let $h^E_\om$ be the metric on $E$ defined by $|1|_{h^E_\om}^2= \eta ^{-1}$,
here $1$ is the canonical unity element of $E$.
Let $\langle\quad \rangle_\om$ be the Hermitian product on
$\cC^\infty (X, L^p\otimes E)= \cC^\infty (X, L^p)$
induced by $h^L, h^E_\om$, $dv_X$ as in \eqref{n2}.
 Then by \eqref{n3},
\begin{equation}\label{n4}
(\cC^\infty (X, L^p\otimes E), \left\langle \quad \right \rangle_\om)
= (\cC^\infty (X, L^p), \left\langle \quad \right \rangle_\nu).
\end{equation}

Observe that $H^0(X, L^p\otimes E)$ does not depend on $g^{TX}$,
$h^L$ or $h^E$.
If $P_{\om,p}(x,x')$, ($x,x'\in X$) denotes the smooth kernel of
the orthogonal projection $P_{\om,p}$ from
$(\cC^\infty (X, L^p\otimes E), \left\langle\,\cdot,\cdot \right \rangle_\om)$
onto $H^0(X, L^p\otimes E)= H^0(X, L^p)$ with respect to $dv_{X}(x)$,
from \eqref{n3}, as in \cite[(3.38)]{MM04a}, we have
\begin{equation}\label{n5}
P_{\nu, p}(x,x')= \eta(x')\, P_{\om,p}(x,x').
\end{equation}
For $f\in \cC^\infty (X)$, set
 \begin{align}\label{n7}
\begin{split}
&K_{\om,p}(x,x')= |P_{\om, p}(x,x')|^2_{(h^{L^p}\otimes h^E_\om)_x\otimes
(h^{L^{p*}}\otimes h^{E^*}_\om)_{x'}},\\
 &(K_{\om,p}f)(x)=  \int_X  K_{\om,p}(x,y)f(y) dv_X(y).
\end{split}
\end{align}

By the definition of the metric $h^E, h^E_\om$,
if we denote by $1^*$ the dual of the section $1$ of $E$,  we know 
\begin{equation}\label{1n8}
1=|1\otimes 1^*|_{h^E\otimes h^{E^*}}^2(x,x')
=|1\otimes 1^*|_{h^E_\om\otimes h^{E^*}_\om}^2(x,x')\eta(x)\eta ^{-1}(x').
\end{equation}
Recall that we identified $(L^p, h^{L^p})$ 
to $(L^p\otimes E,h^{L^p}\otimes h^E)$ by Section $1$. 
Thus from \eqref{1n1}, \eqref{n5} and \eqref{1n8}, we get 
\begin{equation}\label{n8}
K_p(x,x')= |P_{\nu, p}(x,x')|^2
_{(h^{L^p}\otimes h^E)_x\otimes (h^{L^{p*}}\otimes h^{E^*})_{x'}}
= \eta (x)\, \eta (x')\, K_{\om,p}(x,x'),
\end{equation}
and from \eqref{n3}, \eqref{1n2} and \eqref{n8},
\begin{equation}\label{n9}
Q_{K_p} (f)(x)
= \frac{1}{R_p} \int_X  K_{\om,p}(x,y) \eta (x)f(y) dv_X(y).
\end{equation}

Now for the kernel $P_{\om,p}(x,x')$, we can apply
the full asymptotic expansion \cite[Theorem 4.18$^\prime$]{DLM04a}.
In fact let $\overline{\partial} ^{L^p\otimes E,*_\om}$
be the formal adjoint of
the Dolbeault operator $\overline{\partial} ^{L^p\otimes E}$
on the Dolbeault complex
 $\Omega ^{0,\bullet}(X, L^p\otimes E)$ with the scalar product
induced by $g^{TX}$, $h^L$, $h^E_\om$, $dv_X$ as in (\ref{n2}),
and set
\begin{equation}\label{n10}
D_p=\sqrt{2}(\overline{\partial} ^{L^p\otimes E}
+ \overline{\partial} ^{L^p\otimes E,*_\om}).
\end{equation}
Then $H^0(X, L^p\otimes E)= \Ker D_p$ for $p$ large enough,
and $D_p$ is a Dirac operator,
as $g^{TX}(\cdot, \cdot)= \om(\cdot, J\cdot)$ is a K\"ahler metric on $TX$.

Let $\nabla^E$ be the holomorphic Hermitian connection on $(E,
h^E_\om)$. Let $\nabla^{TX}$ be the Levi-Civita connection on $(TX,
g^{TX})$. Let $R^E$, $R^{TX}$ be the corresponding curvatures.

Let $a^X$ be the injectivity radius of $(X, g^{TX})$.
We fix  $\var\in ]0,a^X/4[$. We denote by  $B^{X}(x,\varepsilon)$ and
$B^{T_xX}(0,\varepsilon)$
the open balls in $X$ and  $T_x X$ with center $x$ and radius $\varepsilon$.
We identify $B^{T_xX}(0,\varepsilon)$ with $B^{X}(x,\varepsilon)$
by using the exponential map of $(X, g^{TX})$.

We fix $x_0\in X$.
For $Z\in B^{T_{x_0}X}(0,\var)$ we identify 
$(L_Z, h^L_Z)$, $(E_Z, h^E_Z)$ and $(L^p\otimes E)_Z$
to $(L_{x_0},h^L_{x_0})$, $(E_{x_0},h^E_{x_0})$ and $(L^p\otimes E)_{x_0}$
by parallel transport with respect to the connections
$\nabla ^L$, $\nabla ^E$ and $\nabla^{L^p\otimes E}$ along the curve
$\gamma_Z :[0,1]\ni u \to \exp^X_{x_0} (uZ)$.
Then under our identification, $P_{\om,p}(Z,Z')$ is a function on
$Z,Z'\in T_{x_0}X$, $|Z|,|Z'|\leq \var$,
we denote it by $P_{\om,p,x_0}(Z,Z')$.
Let $\pi : TX\times_{X} TX \to X$ be the
natural projection from the fiberwise product of $TX$ on $X$.
Then
we can view  $P_{\om,p,x_0}(Z,Z')$ as a smooth function on
 $TX\times_{X} TX$ (which is defined for $|Z|,|Z^\prime|\leq \var$)
by identifying a section $S\in \cC^\infty (TX\times_{X}TX,\pi ^* \End (E))$
with the family $(S_x)_{x\in X}$, where
$S_x=S|_{\pi^{-1}(x)}$, $\End(E)=\C$.

We choose $\{ w_i\}_{i=1}^n$
an orthonormal basis of $T^{(1,0)}_{x_0} X$,
then $e_{2j-1}=\tfrac{1}{\sqrt{2}}(w_j+\overline{w}_j)$ and
$e_{2j}=\tfrac{\sqrt{-1}}{\sqrt{2}}(w_j-\overline{w}_j)\,,
 j=1,\dotsc,n\, $
forms an orthonormal basis of $T_{x_0}X$.
We use the coordinates on $T_{x_0}X\simeq\R^{2n}$
where the identification is given by
\begin{equation}\label{n11}
 (Z_1,\cdots, Z_{2n}) \in \R^{2n} \longrightarrow \sum_i
Z_i e_i\in T_{x_0}X.
\end{equation}

In what follows we also introduce the complex coordinates $z=(z_1,\cdots,z_n)$
on $\C^n\simeq\R^{2n}$.
By \cite[(4.114)]{DLM04a} (cf. \cite[(1.91)]{MM04a}), set
\begin{equation}\label{n12}
\begin{split}
P^N(Z,Z') &=\exp\Big(-\frac{\pi}{2}\sum_i
\big(|z_i|^2+|z^{\prime}_i|^2 -2z_i\overline{z}_i^{\prime}\big)\Big).
\end{split}
\end{equation}
Then $P^N$ is the classical Bergman kernel on $\C^n$
(cf. \cite[Remark 1.14]{MM04a}) and
\begin{equation}\label{n13}
|P^N( Z, Z')|^2= e ^{-\pi |Z-Z^\prime|^2}.
\end{equation}

By \cite[Proposition  4.1]{DLM04a},
 for any $l,m\in \N$, $\var>0$, there exists $C_{l,m,\var}>0$
such that for $p\geq 1$, $x,x'\in X$,
\begin{align}\label{1n13}
&|P_{\om,p}(x,x')|_{\cC^m(X\times X)}
\leq C_{l,m,\var}\, p^{-l} \quad \mathrm{if} \,
d(x,x') \geq \var.
\end{align}
Here the $\cC^m$-norm is induced by $\nabla ^L$, $\nabla ^E$, $\nabla ^{TX}$
and $h^L, h^E, g^{TX}$.

By \cite[Theorem 4.18$^\prime$]{DLM04a},
there exist  $J_{r}(Z,Z')$  polynomials
in $Z,Z'$, such that for any $k,m, m'\in \N$,
there exist $N\in \N, C>0, C_0>0$ such that for
$\alpha, \alpha' \in \N^{n}$, $|\alpha|+|\alpha'|\leq m$,
$Z,Z'\in T_{x_0}X$, $|Z|, |Z'|\leq  \var$, $x_0\in X$, $p\geq 1$,
\begin{multline}\label{n14}
\left |\frac{\partial^{|\alpha|+|\alpha'|}}
{\partial Z^{\alpha} {\partial Z'}^{\alpha'}}
\left (\frac{1}{p^n}  P_{\om,p,x_0}(Z,Z')
-\sum_{r=0}^k  (J_{r}P^N) (\sqrt{p} Z,\sqrt{p} Z')
p^{-r/2}\right )\right |_{\cC^{m'}(X)}\\
\leq C  p^{-(k+1-m)/2}  (1+|\sqrt{p} Z|+|\sqrt{p} Z'|)^N
\exp (- C_0 \sqrt{p} |Z-Z'|)+ \cO(p^{-\infty}).
\end{multline}
Here $\cC^{m'}(X)$ is the $\cC^{m'}$ norm for the parameter $x_0\in X$.
The term $\cO(p^{-\infty})$ means that for any $l,l_1\in \N$,
there exists $C_{l,l_1}>0$ such that its $\cC^{l_1}$-norm is dominated
by $C_{l,l_1} p^{-l}$.
(In fact,  by \cite[Theorems 4.6 and 4.17, (4.117)]{DLM04a}
(cf. \cite[Theorem 1.18, (1.31)]{MM04a}), the polynomials $J_{r}(Z,Z')$ have
 the same parity as $r$ and $\deg J_{r}(Z,Z')\leq 3r$,
whose coefficients
are polynomials in $R^{TX}$, $R^E$
and their derivatives of order $\leqslant r-1$).

Now we claim that in \eqref{n14},
\begin{equation}\label{1n14}
J_0=1,\quad J_1(Z,Z')=0.
\end{equation}
In fact, let $dv_{T_{x_0}X}$ be the Riemannian volume form on
$(T_{x_0}X, g^{T_{x_0}X})$, and $\kappa$ be the function defined by
\begin{equation}\label{1n15}
dv_X(Z)= \kappa(x_0,Z) dv_{T_{x_0}X}(Z).
\end{equation}
Then (also cf. \cite[(1.31)]{MM04a})
\begin{equation}\label{1n16}
\kappa(x_0,Z)= 1
+ \frac{1}{6} \left\langle R^{TX}_{x_0} (Z,e_i) Z, e_i\right\rangle_{x_0}
 + \cO (|Z|^3).
\end{equation}

As we only work on $\cC^\infty(X, L^p\otimes E)$,
by \cite[(4.115)]{DLM04a}, we get the first equation in \eqref{1n14}.

Recall that in the normal coordinate, after the rescaling $Z\to Z/t$
with $t=\frac{1}{\sqrt{p}}$, we get an operator $\cL_t$
from the restriction of $D_p^2$ on $\cC^\infty(X, L^p\otimes E)$
which has the following formal expansion
 (cf. \cite[(1.104)]{DLM04a}, \cite[Theorem 1.4]{MM04a}),
\begin{equation}\label{1n17}
\cL_t= \cL+ \sum_{r=1}^\infty \mathcal{Q}_r t^r.
\end{equation}
Now, from \cite[Theorem 5.1]{DLM04a} (or \cite[(1.87), (1.97)]{MM04a}),
\begin{align}\label{1n18}
\cL= \sum_{j=1}^n (-2{\tfrac{\partial}{\partial z
_i}}+\pi \overline{z}_i)
(2{\tfrac{\partial}{\partial\overline{z}_i}}+\pi z_i),\quad \mathcal{Q}_1=0.
\end{align}
(In fact, $P^N(Z,Z^\prime)$ is the smooth kernel of the orthogonal
projection from $L^2(\R)$ onto $\Ker (\cL)$).
Thus from \cite[(4.107)]{DLM04a} (cf. \cite[(1.111)]{MM04a}),
\eqref{1n16} and \eqref{1n18} we get the second equation of \eqref{1n14}.

Note that
$|P_{\om,p,x_0}(Z,Z')|^2 =P_{\om,p,x_0}(Z,Z')\overline{P_{\om,p,x_0}(Z,Z')}$,
thus from \eqref{n7}, \eqref{n14} and \eqref{1n14},
there exist  $J^\prime_{r}(Z,Z')$  polynomials
in $Z,Z'$ such that
\begin{multline}\label{n15}
\left |\frac{1}{p^{2n+1}} \Delta_Z
\Big(K_{\om,p,x_0}(Z,Z')
- \Big(1+  \sum_{r=2}^k p^{-r/2}  J^\prime_{r} (\sqrt{p} Z,\sqrt{p} Z')\Big)
e ^{-\pi p |Z-Z^\prime|^2} \Big)\right|\\
\leq C p^{-(k+1)/2}
(1+|\sqrt{p} Z|+|\sqrt{p} Z'|)^N
\exp (- C_0 \sqrt{p} |Z-Z'|)+ \cO(p^{-\infty}).
\end{multline}

For a function $f\in \cC^\infty(X)$, we denote it as $f(x_0,Z)$
a family (with parameter $x_0$) of function on $Z$ in the normal coordinate
near $x_0$.
Now, for any polynomial $Q_{x_0}(Z^\prime)$,
we define the operator
\begin{align}\label{n16}
(\mQ_p f)(x_0)
=  p^n \int_{|Z^\prime|\leq \var}
Q_{x_0}(\sqrt{p}Z^\prime) e^{-\pi p |Z^\prime|^2} f(x_0,Z^\prime)
dv_X(x_0,Z^\prime).
\end{align}
Then we observe that there exists $C_1>0$ such that for any
$p\in \N, f\in \cC^\infty(X)$, we have
\begin{align}\label{n17}
|\mQ_p f|_{L^2} \leq C_1 |f|_{L^2}.
\end{align}
In fact,
\begin{multline}\label{n18}
|\mQ_p f|_{L^2}^2
\leq \int_X dv_X(x_0) \Big\{
p^n \Big (\int_{|Z^\prime|\leq \var}
|Q_{x_0}(\sqrt{p}Z^\prime)| e^{-\pi p |Z^\prime|^2}dv_X(x_0,Z^\prime)\Big)\\
\times p^n\Big (\int_{|Z^\prime|\leq \var}
|Q_{x_0}(\sqrt{p}Z^\prime)| e^{-\pi p |Z^\prime|^2}
|f(x_0,Z^\prime)|^2
dv_X(x_0,Z^\prime)\Big)\Big\}\\
\leq C' \int_X dv_X(x_0) p^n\int_{|Z^\prime|\leq \var}
|Q_{x_0}(\sqrt{p}Z^\prime)| e^{-\pi p |Z^\prime|^2}
|f(x_0,Z^\prime)|^2
dv_X(x_0,Z^\prime)\\
\leq C_1 |f|_{L^2}^2.
\end{multline}

Observe that in the normal coordinate, at $Z=0$,
$\Delta_Z= -\sum_{j=1}^{2n} \tfrac{\partial^2} {\partial Z_j^2}$. Thus
\begin{equation} \label{n19}
(\Delta_Z e^{-\pi p |Z-Z^\prime|^2})|_{Z=0}
= 4\pi p(n- \pi p |Z^\prime|^2) e^{-\pi p |Z^\prime|^2}.
\end{equation}
Thus from \eqref{n13}, \eqref{n14}, \eqref{1n14}, \eqref{n15} and \eqref{n17},
we get
\begin{align}\label{n21}
\begin{split}
&\left|p^{-n}  K_{\om,p} f
- p^n \int_{|Z^\prime|\leq \var}e^{-\pi p |Z^\prime|^2}f(x_0,Z^\prime)
dv_X(x_0,Z^\prime)\right|_{L^2}
\leq \frac{C}{p}  \left|f\right|_{L^2},\\
&\left|p^{-n-1}\Delta  K_{\om,p} f
- 4\pi p^n \int_{|Z^\prime|\leq \var}
(n- \pi p |Z^\prime|^2) e^{-\pi p |Z^\prime|^2}f(x_0,Z^\prime)
dv_X(x_0,Z^\prime)\right|_{L^2}
\leq \frac{C}{p} \left|f\right|_{L^2}.
\end{split}
\end{align}

Set 
\begin{align}\label{1n21}
\begin{split}
&K_{\eta, \om,p}(x,y)
= \langle d\eta (x), d_x K_{\om,p}(x,y) \rangle_{g^{T^*X}},\\
&(K_{\eta, \om,p}f)(x) = \int_X K_{\eta, \om,p}(x,y) f(y) dv_X(y).
\end{split}
\end{align}
Then from \eqref{n14}, \eqref{1n14} and \eqref{n17}, we get
\begin{align}\label{1n23}
\left|p^{-n-1} K_{\eta, \om,p}f -2\pi p^n  \int_{|Z^\prime|\leq \var}
\sum_{i=1}^{2n} (\frac{\partial}{\partial Z_i}\eta)(x_0,0) Z_i^\prime
e^{-\pi p |Z^\prime|^2}f(x_0,Z^\prime)
dv_X(x_0,Z^\prime)\right|_{L^2}\leq \frac{C}{p} \left|f\right|_{L^2}.
\end{align}

Let $e ^{-u\Delta}(x,x')$ be the smooth kernel of the heat operator
$e^{-u \Delta}$ with respect to $dv_X(x')$.
Let $d(x,y)$ be the Riemannian distance from $x$ to $y$ on $(X, g^{TX})$.
By the heat kernel expansion in
\cite[Theorems 2.23, 2.26]{BeGeVe}, there exist $\Phi_i(x,y)$
smooth functions on $X\times X$ such that when $u\to 0$,
we have the following asymptotic expansion 
\begin{align}\label{1n22}
\left| \frac{\partial ^l}{\partial u^l}
\Big(e^{-u\Delta}(x,y) -(4\pi u)^{-n} \sum_{i=0}^k u^i\Phi_i(x,y)
e^{-\frac{1}{4u}d(x,y)^2}\Big)\right|_{\cC^m(X\times X)} 
= \cO(u^{k-n-l-\frac{m}{2}+1}),
\end{align}
and
\begin{align}\label{n23}
\Phi_{0}(x,y)=1.
\end{align}

If we still use the normal coordinate, then by \eqref{1n22},
there exist $\phi_{i,x_0}(Z^\prime):=\Phi_i(0,Z^\prime)$ such that uniformly
for $x_0\in X$, $Z^\prime\in T_{x_0}X, |Z^\prime|\leq \var$,
we have the following asymptotic expansion when $u\to 0$,
\begin{align}\label{n22}
\left| \frac{\partial ^l}{\partial u^l} \Big(e^{-u\Delta}(0, Z')
- (4\pi u)^{-n} \Big(1+ \sum_{i=1}^k u^i \phi_{i,x_0}(Z')\Big)
e^{-\frac{1}{4u} |Z^\prime|^2}\Big)\right| =\cO(u^{k-n-l+1}),
\end{align}
and 
\begin{multline}\label{1n24}
\Big|\langle d\eta (x_0), d_{x_0}e^{-u\Delta} 
\rangle_{g^{T^*X}}(0,Z^\prime)\\
-  (4\pi u)^{-n} 
\sum_{i=1}^{2n} (\frac{\partial}{\partial Z_i}\eta)(x_0,0)\frac{Z_i^\prime}{2u}
 \Big(1+ \sum_{i=1}^k u^i \phi_{i,x_0}(Z')\Big)\Big) 
e^{-\frac{1}{4u} |Z^\prime|^2}\\
-  (4\pi u)^{-n} \sum_{i=1}^k u^i
\langle d\eta (x_0), (d_{x_0}\Phi_i)(0,Z^\prime)\rangle
 e^{-\frac{1}{4u} |Z^\prime|^2}
 \Big| = \cO(u^{k-n+\frac{1}{2}}).
\end{multline}

Observe that
\begin{align}\label{n24}
\frac{1}{p} \Delta \exp\Big(-\frac{\Delta}{4\pi p}\Big)
= - \frac{1}{p} (\tfrac{\partial}{\partial u} e^{-u\Delta})
|_{u= \frac{1}{4\pi p}}.
\end{align}

Now from \eqref{n17}, \eqref{n21}--\eqref{n24},  we get
\begin{align}\label{n25}
\begin{split}
&\left|\Big(p^{-n}  K_{\om,p}
- \exp\Big(-\frac{\Delta}{4\pi p}\Big)\Big) f \right|_{L^2}
\leq \frac{C}{p}  \left|f\right|_{L^2},\\
&\left|\frac{1}{p} \Big(p^{-n}\Delta  K_{\om,p} -
\Delta\exp\Big(-\frac{\Delta}{4\pi p}\Big)\Big) f \right|_{L^2}
\leq \frac{C}{p}  \left|f\right|_{L^2}.
\end{split}
\end{align}
and
\begin{align}\label{1n25}
&\left|\frac{1}{p}\Big(p^{-n}K_{\eta,\om,p} 
- \langle d\eta, d \exp({-\frac{\Delta}{4\pi p}}) \rangle \Big)f \right|_{L^2}
\leq \frac{C}{p}  \left|f\right|_{L^2}.
\end{align}

Note that
\begin{multline}\label{n26}
(\Delta \eta K_{\om,p})(x,y)
= (\Delta \eta)(x)K_{\om,p}(x,y) +\eta(x) \Delta_x K_{\om,p}(x,y)\\
-2 \langle d\eta (x), d_x K_{\om,p}(x,x') \rangle_{g^{T^*X}},
\end{multline}
and 
$R_p= \frac{{\rm Vol}(X, dv_X)}{{\rm Vol}(X, \nu)} p^n + \cO(p^{n-1})$.
 From \eqref{n9}, \eqref{n25}-\eqref{n26}, we get \eqref{1n4}.

To get the last part of Theorem \ref{nt1}, as we noticed in \cite[\S
4.5]{DLM04a}, the constants in \eqref{n14} will be uniformly bounded
under our condition, thus we can take $C$ in \eqref{1n4}, 
\eqref{n25}and  \eqref{1n25}
independent of $h^L$, $d\nu$.
\end{proof}

\subsection*{Acknowledgments}
We thank Professor Simon Donaldson for useful communications.

\def\cprime{$'$} \def\cprime{$'$}
\providecommand{\bysame}{\leavevmode\hbox to3em{\hrulefill}\thinspace}
\providecommand{\MR}{\relax\ifhmode\unskip\space\fi MR }
\providecommand{\MRhref}[2]{%
  \href{http://www.ams.org/mathscinet-getitem?mr=#1}{#2}
}
\providecommand{\href}[2]{#2}


\begin{thebibliography}{10}


\bibitem{BeGeVe}
N.~Berline, E.~Getzler, and M.~Vergne, \emph{Heat kernels and {D}irac
  operators}, Grundl. Math. Wiss. Band 298, Springer-Verlag, Berlin, 1992.

\bibitem{DLM04a}
X.~Dai, K.~Liu, and X.~Ma, \emph{On the asymptotic expansion of {B}ergman
  kernel}, C. R. Math. Acad. Sci. Paris \textbf{339} (2004), no.~3, 193--198.
The full version: J. Differential Geom. to appear, math.DG/0404494.

\bibitem{D05}
S.~K. Donaldson,
\emph{Some numerical results in complex differential geometry},
math.DG/0512625.

\bibitem{MM04a}
X.~Ma and G.~Marinescu,
\emph{Generalized {B}ergman kernels on symplectic manifolds},
C. R. Math. Acad. Sci. Paris \textbf{339} (2004), no.~7, 493--498.
The full version: math.DG/0411559.

\bibitem{MM06a}
\bysame, \emph{Holomorphic Morse Inequalities and Bergman Kernels},
book in preparation, (2006).


\end{thebibliography}
\end{document}